\documentclass[12]{report}
\usepackage{graphicx}
\usepackage{cancel}
\usepackage{amsmath}    

\date{}

\begin{document}
\title{An Approximation for Square Roots known in India prior to Alexander's Conquest\footnote{Translated by Dileep Karanth, Department of Physics, University of Wisconsin-Parkside, Kenosha.}}
\author{L\'eon Rodet}
\maketitle

\long\def\symbolfootnote[#1]#2{\begingroup%
\def\thefootnote{\fnsymbol{footnote}}\footnote[#1]{#2}\endgroup} 

\begin{abstract}
This paper was published in French under the title  \textit{Sur une m\'ethode d'approximation des racines carr\'ees dans l'Inde ant\'erieurement \`a la conqu\^ete d'Alexandre}, in the  Bulletin de la S. M. F., volume 7 (1879), p. 98-102. 

It argues that the process of approximating the irrational square root of an integer, now commonly known as Newton's method, was known in India since very ancient times.

\end{abstract}

In a Treatise on Arithmetic that he composed in 613 AH or 1216 CE ``in response to the requests of his friends", the Persian author Hassan ben al-Hossein al-Hakak al-Morouzi (that is to say, native of Merv in Khorasan), gives the following procedure to obtain an approximation to an irrational square root. First one obtains an approximation to the desired square root to within an integer from below, by means of the well-known method, then ``one divides the remainder of the operation by the double of the root obtained \textit{plus one}, and the root is composed of the integral part and this fraction". Thus
\begin{equation*}
\sqrt{12} = 3 \frac{3}{7}  \mbox{,    } \sqrt{145} = 12 \frac{1}{25}
\end{equation*}

This procedure is not hard to substantiate. Let us we denote by $N$ the given number, by $a$ its root to the nearest integer, and by $r$ the remainder, so that $ N = a^2 + r$; Let, in addition, $\epsilon$ be complement of the root, so that
 
\begin{equation*}
r = (2a + \epsilon)\epsilon  \mbox{,       or         } \epsilon = \frac{r}{2a + \epsilon};
\end{equation*}
 Since $\epsilon  < 1$ always, one can approximate from below by replacing it by one in the denominator of the fraction, that is, by setting
\begin{equation*}
\epsilon_1   = \frac{r}{2a + 1};
\end{equation*}

	This is the first time that I have found this approximation method in the writings of an oriental author, and though I should have found it quite remarkable for for the time when he lived, I would probably not have given it much attention if I had not been led to think that this procedure must have been practised quite commonly in India, even since a very early period. 
	
	Since Al-Morouzi has given the approximation $ \sqrt{12} = 3 \frac{3}{7}$, I was naturally led to see whether in this system, $ \sqrt{10} = 3 \frac{1}{7} = \frac{22}{7}$, that is, the value given by Archimedes for the number $\pi$ -- the sole value in usage for a long time among the Greeks. Indeed, it is the only one used by Hero the Younger in his Treatise of surfaces and volumes. 
	
	However, since $ \frac{22}{7} = \sqrt{10}$, we now understand why the Indian astronomer Brahmagupta used to teach his disciples (see Colebrooke's translation, number 40, second part of the rule) that ``the root of 10 times the half-diameter (radius), multiplied by this half-diameter itself, is the area of the circle." In the first part of this rule, he multiplied the diameter in the first case, and the square of the radius in the second case, by 3, and he regarded this process as ``usual", or ``popular" (\textit{vy\^avah\^arika}); On the other hand, the second, which sets $\pi = \sqrt{10}$, is called ``\textit{s\^uxma}", ``fine, subtle" and even possibly ``refined". The question has been asked, by myself after many others before me, where Brahmagupta had found this value for $\pi$, which he regarded as exact, so to say. The Persian writer whom I have cited has given us an answer to this question. It is simply that $ \pi= \frac{22}{7}$ is an approximation \textit{from above}; $ \sqrt{10} = \frac{22}{7}$ is an approximation \textit{from below}. Perhaps Brahmagupta, in his oral teaching, would have known to make this distinction to his students.
		
		As the example borrowed from Brahmagupta gave me to think that the Indians had known and practised the method of approximation of square roots that I have mentioned, the idea came to me to examine from this point of view a value of $\sqrt{2}$ given in the \textit{String Precepts}\footnote{Translator's Note: This work is more commonly referred to by its untranslated title: \textit{\'Sulba S\=utra.}} of Baudh\^ayana, a curious work composed probably in the IVth century before our era, and which contains the rules that Brahmins have to observe in constructing their altars.  This value of  $\sqrt{2}$, with whose interpretation I was equally preoccupied, is

\begin{equation*}
\sqrt{2} = 1 + \frac{1}{3} + \frac{1}{3 . 4} -  \frac{1}{3 . 4 . 34}
\end{equation*}
This form of series seemed to me to be completely original. 

Now, if we apply the method in question to $\sqrt{2}$, we would have to set 
\begin{equation*}
N =2 \mbox{,      } a =1\mbox{,      }r = 1\mbox{,      } \epsilon_1 =  \frac{1}{2 + 1 }=  \frac{1}{3}
\end{equation*}
Thus the first fractional term in the series given by our Brahmin is thus obtained by the method of the Persian; $1 + 1 \frac{1}{3}$ is a an approximation \textit{from below} for the desired root. 
	Let us continue the operation. We have to calculate
\begin{equation*}
r^{\prime}  = r - (2a + \epsilon_1)\epsilon_1 \mbox{,      } \epsilon_2 = \frac{r^{\prime}}{(2a + \epsilon_1)+ \epsilon_2}
\end{equation*}
which gives us
\begin{equation*}
r^{\prime}  = 1 - (2 + \frac{1}{3})\frac{1}{3} = 1 - \frac{7}{9} = \frac{2}{9},
\end{equation*}
\begin{equation*}
\epsilon_2 = \frac{2}{9}: [ 2 (1 + \frac{1}{3}) + \epsilon_2 ]   =  \frac{2}{9} : (\frac{8}{3} + \epsilon_2).
\end{equation*}

Now our second term is 
\begin{equation*}
\frac{1}{3.4} = \frac{2}{9} . \frac{3}{8}
\end{equation*}
and this time, the author has neglected $  \epsilon_2 $ in the denominator of the divisor fraction and obtained an approximation for $ \sqrt{2}$ \textit{from above} in the series

\begin{equation*}
1 + \frac{1}{3} + \frac{1}{3 . 4} 
\end{equation*}

If we pursued this operation further, we would find that the quantity to be subtracted from $ r^{\prime} = \frac{2}{9}$ to find $r^{\prime \prime}$ is $ \frac{33}{144}$. Now $\frac{2}{9} =  \frac{33}{144}$; hence $r^{\prime \prime} = - \frac{1}{144}$. 
I give this result in this negative form, because I have already observed elsewhere that the Indians developed the notion of \textit{negative numbers} very early, and did not concern themselves with the sign of the result, in their calculations.

On the other hand, the double of the new root being $ \frac{33}{12}$, we will have
\begin{equation*}
\epsilon_3 = (- \frac{1}{144})/\frac{34}{12} = -\frac{1}{3.4.34}.
\end{equation*}
Note that, this corrective term being too high in absolute value, the expression 
\begin{equation*}
1 + \frac{1}{3} + \frac{1}{3 . 4} -  \frac{1}{3 . 4 . 34}
\end{equation*}
is again approximated \textit{from below}.

Nothing prevents the continuation of this calculation, and there is no indication that the authors of this series of four terms could not have pursued the operation, if they had any need to. 

From this we can conclude that, at the time that Baudh\^ayana lived, that is around the IVth century before our era, people knew of an expression for the value of a surd square root in the form of a series approximating the root alternatingly from from below and from above. The method of obtaining the first term of this series was as follows: First the square root of the number is found, to the closest integer. Its square is subtracted from the given number. The remainder is then divided by twice the root plus one. The following terms, alternatingly positive and negative, are obtained by a process none other than that which is known today by the name 
\textit{Newton's Approximation Method}. 

In conclusion, I would add that Al-Mourouzi applied the same process (limited to the first term in the series) to cube roots. Here $ \epsilon = \frac{r}{3a^2 + 3a \epsilon + \epsilon^2} $; again he sets $ \epsilon = 1$ and gets $\epsilon_1 = \frac{r}{3a^2 + 3a  + 1}$. As a numerical example, he gives 
 \begin{equation*}
\sqrt[3]{1748} = 12 \frac{20}{469}
\end{equation*}
However, I am not aware if one could find examples of the application of this process in the writings of Indian authors, or whether for the square root, they continued to use a series expansion, and whether the following terms of the series were of the form $\epsilon_n = \frac{r_n}{3a^3_{n-1}} $, as in the method of Newton.

\section*{Translator's Acknowledgements}
The translator is grateful to Dr. Albrecht Heeffer for introducing him to L\'eon Rodet's writings, and for encouraging him to translate them. The original French article was made available by www.numdam.org, and is freely available for download.   

\end{document}